\documentclass[12pt]{amsart}

\usepackage{latexsym}
\usepackage{psfrag}
\usepackage{amsmath}
\usepackage{amssymb}
\usepackage{epsfig}
\usepackage{amsfonts}

\newcommand{\excise}[1]{}

\newtheorem{thm}{Theorem}[section]

\newtheorem{conv}[thm]{Convention}

\newtheorem{Warn}[thm]{Caution}

    {\end{Example}}
    {\end{Remark}}

\def\sq{\square}

\def\rr{\mathbb R}

\def\sm{\smallsetminus}

\def\Ga{\Gamma}

\def\ssu{\subset}

\def\<{\langle}
\def\>{\rangle}

\def\0{{\mathbf 0}}

\def\.{\hskip.06cm}
\def\ts{\hskip.03cm}

\def\pt{\partial}

\def\bio{{\text{\bf {0}}}}

\def\bbe{e}
\def\bw{\textbf{\textit{w}}}

\def\bv{\textbf{\textit{v}}}
\def\ba{\textbf{\textit{a}}}

\begin{document}

\title[Short proof of the infinitesimal rigidity]%
   {A short proof of rigidity of convex
polytopes}
\author[Igor~Pak]{Igor~Pak}
\date{11 August 2005}



\begin{abstract}
\noindent
We present a much simplified proof of Dehn's theorem
on the infinitesimal rigidity of convex polytopes.  Our
approach is based on the ideas of
Trush-kina~\cite{Tru} and Schramm~\cite{Sch}.
\end{abstract}

\maketitle

\vspace{-1ex}

\part*{Introduction}\label{intropage}

\noindent
Let~$P \ssu \rr^3$ be a simplicial convex polytope.
Define a \emph{continuous deformation} $\{P_t : t \in [0,1]\}$
of~$P = P_0$ to be a family of convex polytopes with the same
combinatorial structure, the same corresponding edge lengths,
and continuity on vertices.  We say that~$P$ is
\emph{continuously rigid} if every such deformation
is a rigid motion in~$\rr^3$.  It is a classical corollary
from the Cauchy theorem that every simplicial polytope~$P$
is continuously rigid (see below).  In this article we
present a simple proof of Dehn's theorem, which also
implies the continuous rigidity.

\smallskip

Let $\bv_i(t) = \.\stackrel{\longrightarrow}{Ov_i}$ be a vector
from the origin~$O$ into the vertex~$v_i$ in~$P_t$.  Think of
vectors~$\bv_i'(t)$ as of \emph{velocities} of
vertices~$v_i$.  For an edge length $|v_iv_j|$ to be constant
under the deformation we need $\|\bv_i(t)-\bv_j(t)\|' =  0$,
where $\|\bw\| = (\bw,\bw) = |\bw|^2$.
Thus, in particular, at $t=0$ we have:
$$\aligned
0\. & = \frac{d}{d\ts t}\,\bigl\|\bv_i(t)-\bv_j(t)\bigr\|_{t=0}
\.  = \. \frac{d}{d\ts t}\,
\bigl\|\bigl(\bv_i(0)-\bv_j(0)\bigr) +
t\ts \bigl(\bv_i'(0)-\bv_j'(0)\bigr)\bigr\|_{t=0}\\
\. & = \. 2\.\bigl(\bv_i(0)-\bv_j(0),\bv_i'(0)-\bv_j'(0)\bigr).
\endaligned
$$
This leads to the following definition of infinitesimal rigidity.

Let $V= \{v_1,\ldots,v_n\}$ and~$E$ be the set of vertices
and edges of~$P$.  We will assume that vertices $v_1,v_2$
and~$v_3$ form a triangular face of~$P$; they are called
\emph{base vertices}, and face~$(v_1v_2v_3)$ is called
\emph{base triangle}.  Suppose we are given a
vector~$\ba_i \in \rr^3$, for every $i \in \{1,\ldots,n\}$.
We say that the set of vectors $\{\ba_i, 1 \le i \le n\}$ defines an
\emph{infinitesimal rigid motion} if
$$
(\ast) \ \qquad
(\bv_i-\bv_j,\ba_i-\ba_j) \. = \. 0, \ \ \, \text{for every}
\ \, (v_i,v_j) \in E.
$$
An infinitesimal rigid motion is called \emph{planted}
if the base velocities vectors are equal to zero:
$\ba_1 = \ba_2 = \ba_3 = \bio$.   Finally, we say
that a simplicial polytope~$P \ssu \rr^3$ is
\emph{infinitesimally rigid} if every planted
infinitesimal rigid motion is trivial: $\ba_i = \bio$,
for all~$i \in \{1,\ldots,n\}$.  Here is the main result
of this paper:

\medskip

\noindent
{\bf Dehn's Theorem.} \ {\it
Every simplicial convex polytope in~$\ts \rr^3$ is infinitesimally
rigid.}

\bigskip

Of course, the restriction to planted infinitesimal
rigid motions is necessary, as the usual rigid motions
of~$P$ in~$\rr^3$ can define nontrivial infinitesimal
rigid motions.  By the argument above, the infinitesimal
rigidity implies the continuous rigidity:

\medskip

\noindent
{\bf Corollary} (Cauchy). \ {\it
Every simplicial convex polytope in~$\rr^3$ is
continuously rigid.}

\medskip

In this paper we present a new proof of Dehn's theorem,
based on the approach by Trushkina~\cite{Tru} (see also a
followup~\cite{Tru-next}).  Unfortunately, the technical
details in Trushkina's paper are somewhat complicated.
We substitute Trushkina's definition of the
\emph{inversion} with the one given by Schramm~\cite{Sch}
(for different purposes).
We should mention an important survey~\cite{C3} (see~$\S 4.6$),
which helped us to translate the ideas in~\cite{Sch}
into the language of infinitesimal rigidity.

Before we conclude, let us repeat that the continuous
rigidity of convex polytopes follows from the
classical Cauchy theorem~\cite{Al-book,Ber,C3}.  Dehn's theorem
was established by Dehn in~\cite{Dehn} in an equivalent
language of the \emph{static rigidity}.
This result became fundamental in the modern study of
rigidity of frameworks and non-convex
polyhedra~\cite{C3,Whi-survey}, and a number of
proofs have been found.  We refer to~\cite{Al-inverse,Pak}
for the exposition of Dehn's original proof,
to~\cite{FP,G,Roth,Whi} for applications and
modern treatment, and to~\cite{Al-inverse,C3}
for further references.

\bigskip

\section{Proof of infinitesimal rigidity}

First, note that equations~$(\ast)$ above say that
the difference in velocities along an edge is orthogonal
to this edge of the polytope.  Think of velocity
vectors as vector functions on vertices of~$P$ which
are equal to~$\bio$ on base vertices~$v_1,v_2,v_3$.
The idea of the proof is to enlarge the set of such
functions and prove a stronger result.

As before, let~$V = \{v_1,\ldots,v_n\}$ and~$E$ be the set
of vertices and edges of a simplicial convex polytope
$P\ssu \rr^3$.  Consider the set of all vector sequences
$(\ba_1,\ldots,\ba_n)$, $\ba_i \in \rr^3$, such that
for every edge $(v_i,v_j) \in E$ we have
one of the following three possibilities:

\smallskip

{\footnotesize \bf 1.} \ $(\bv_i-\bv_j,\ba_i) = (\bv_i-\bv_j,\ba_j) = 0$,

{\footnotesize \bf 2.}  \ $(\bv_i-\bv_j,\ba_i) <0$ \. and \.
$(\bv_i-\bv_j,\ba_j) < 0$,

{\footnotesize \bf 3.}  \ $(\bv_i-\bv_j,\ba_i)> 0$ \. and \.
$(\bv_i-\bv_j,\ba_j) > 0$.

\smallskip

\noindent
In other words, we require that projections of velocity
vectors~$\ba_i$ and~$\ba_j$ onto edge $(v_i,v_j)$
have the same signs.
We say that a vertex~$v_i$ is \emph{dead} if
$\ba_i = \bio$; it is \emph{live} otherwise.
We need to prove that for every vector sequence
$(\ba_1,\ldots,\ba_n)$ as above, if the base
vertices are dead, then all vertices~$v_i \in V$ are dead.
By definition of the infinitesimal rigidity,
this would immediately imply the theorem.

Denote by~$\Ga = (V,E)$ the graph of~$P$.  Since~$P$
is simplicial, $\Ga$ is a plane triangulation.
Consider an orientation of edges of~$\Ga$ in the direction
of projections of the velocity vectors. More precisely,
we orient edges $v_i \to v_j$ in case~{\footnotesize \bf 2},
orient them $v_i \gets v_j$ in case~{\footnotesize \bf 3},
and leave them unoriented in case~{\footnotesize \bf 1}.
Clearly, the edges adjacent to dead vertices
are unoriented.

Consider two edges $e=(v_i,v_j)$ and $e'=(v_i,v_r)$,
$e,e' \in E$, with a common vertex~$v_i$,
such that $(v_iv_jv_r)$ is a face in~$P$.
We say that edges~$e$ and~$e'$

\smallskip

$\circ$ \  have \emph{one inversion} \. if one of them
is oriented into~$v_i$, and the other out of~$v_i$,

$\circ$ \  have \emph{zero inversions} \. if both of them
are oriented into~$v_i$ or out of~$v_i$,


$\circ$ \ have a \emph{half-inversion} \. if one of the edges is
oriented and the other is unoriented,


$\circ$ \ have \emph{one inversion} \. if both of them
are unoriented and $v_i$ is a live vertex,

$\circ$ \ have \emph{zero inversions} \. if both of them
are unoriented and $v_i$ is a dead vertex.

\smallskip

\noindent
When we talk about the number of inversions in a subgraph,
in a triangle, or around a vertex, we mean the total
sum of inversions between pairs of edges involved.
For example, we say that a graph \emph{has at least~$q$
inversions} if this sum is~$\ge q$.

\begin{figure}[hbt]
\psfrag{1}{$v_i$}
\psfrag{2}{$v_j$}
\psfrag{3}{$v_r$}
\begin{center}
\epsfig{file=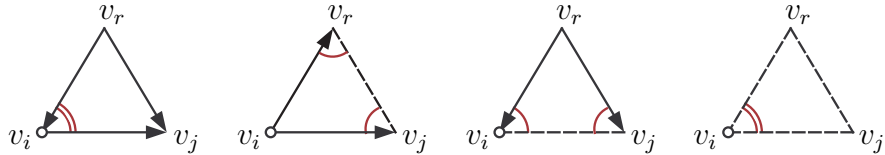,width=11.8cm}
\end{center}
\caption{Different orientations of $(v_iv_jv_r)$, where
vertex~$v_i$ is live.}
\label{f:alg-dehn-schramm}
\end{figure}

A triangle is called \emph{active} if at least one of
its vertices is live; it is called \emph{inactive} otherwise.
Now consider different orientations of an active
triangle $(v_iv_jv_r)$ where vertex~$v_i$ is live
(see some of them in Figure~\ref{f:alg-dehn-schramm}).
A simple enumeration of all possible cases gives the
following result:

\medskip

\noindent
{\bf Lemma 1.} \
{\it Every active triangle has at least one inversion.}


%

This lemma gives a lower bound on the number of inversions
in~$\Ga$.  To get an upper bound, we count inversions
around vertices:

\medskip

\noindent
{\bf Lemma 2.} \
{\it There are at most two inversions around every live vertex.}

\medskip

We postpone the proof of the lemma until we finish the proof
of the theorem.
Consider what this gives us when the only unoriented edges
in~$\Ga$ are the base edges $(v_1,v_2)$, $(v_1,v_3)$, and~$(v_2,v_3)$.
In this case we have~$(n-3)$ live vertices.  Recall also that there
are no inversions around dead vertices.  Thus,
by Lemma~2, there are at most $2(n-3) = 2n-6$ inversions in~$\Ga$.
On the other hand, recall that a triangulation with~$n$ vertices
always has~$(2n-4)$ triangles (see e.g.~\cite{Ber}), and, by
assumption, only one of them is inactive.  Thus, by
Lemma~1, there are at least~$(2n-5)$ inversions in~$\Ga$,
a contradiction.

We use the same strategy in general case.  Remove
from~$\Ga$ all inactive triangles together with all edges
and vertices which belong only to inactive triangles.
Denote by~$H=(V',E')$ a connected component of the
remaining graph.  Since~$\Ga$ is planar, the induced
subgraph~$H$ of~$\Ga$ has a well defined boundary~$\pt H$.
Denote by~$k$ the number of vertices in~$\pt H$
(all of them dead), and by~$\ell$ the number of
connected components of~$\pt H$.  Finally, denote by~$m$
the number of vertices in~$H \sm \pt H$ (some of them
live and some possibly dead).
By Lemma~2, there are
at most~$2\ts m$ inversions in~$H$.

Let us now estimate the number of inversions via the
number~$t$ of triangles in~$H$.  Observe that the total
number of vertices and edges in~$H$ is
given by
$$|V'| \. = \. m+k, \quad \text{and} \ \ \
2\. |E'| \. = \. k+3\ts t\..$$
On the other hand, there are $f = (t+\ell)$
faces in~$H$ taken together with the boundary components.
Substituting these values into Euler's formula $|V'|-|E'|+ f = 2$,
we conclude that graph~$H$ has exactly
$t = 2\ts m + k + 2\ts\ell- 4$ triangles.
Since there is at least one inactive triangle $(v_1v_2v_3)$,
we have~$k \ge 3$ and~$\ell \ge 1$.  Therefore, by
Lemma~1, there are at least
$$t \. = \. 2\ts m \. + \. k \. + \. 2\ts\ell \. - \. 4 \.
\ge \. 2\ts m \. + \. 3 \. + \. 2 \. - \. 4 \. = \.
 2\ts m \. + \. 1
$$
inversions in~$H$, a contradiction. \ $\sq$

\medskip

\section{Proof of Lemma~2.}

Let us consider all possibilities one by one, and check the
claim in each case. Suppose a vertex~$v_i$ is adjacent to
three or more unoriented edges.  This means that~$\ba_i$
is orthogonal to at least three vectors spanning~$\rr^3$.
Therefore, $\ba_i = \bio$ and~$v_i$ is a dead vertex
with zero inversions.

\begin{figure}[hbt]
\begin{center}
\psfrag{0}{\footnotesize \bf 0}
\psfrag{1}{\footnotesize \bf 1}
\psfrag{2}{\footnotesize \bf 2}
\epsfig{file=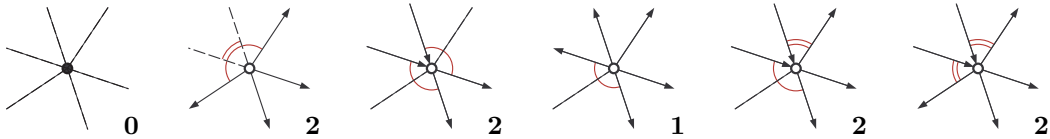,width=14.2cm}
\end{center}
\caption{The number of inversions around a vertex in
different cases.}
\label{f:alg-dehn-inversions}
\end{figure}

Suppose now that~$v_i$ is adjacent to exactly two unoriented
edges~$\bbe,\bbe'$ in~$P$.  This means that~$\ba_i \ne \bio$ is
orthogonal to a plane spanned by these edges.  Observe
that~$\bbe,\bbe'$ separate the edges oriented into~$v_i$
from those oriented out of~$v_i$.  Thus, there are either
two half-inversions and one inversion if the edges~$\bbe,\bbe'$
are adjacent, or four half-inversions if~$\bbe,\bbe'$ are
not adjacent (see Figure~\ref{f:alg-dehn-inversions}).

Next, suppose that~$v_i$ is adjacent to exactly one unoriented
edge~$\bbe$.  Since $\ba_i \ne \bio$ in this case,
consider a plane containing~$v_i$ and orthogonal to~$\ba_i$.
This plane contains~$\bbe$ and separates the edges
oriented into~$v_i$ from those oriented out of~$v_i$.
Therefore, there are either two half-inversions
if all other edges are oriented into~$v_i$, or out
of~$v_i$,  or two half-inversions and one inversion
otherwise.

Finally, if~$v_i$ is a live vertex that is not adjacent
to any unoriented edges, then the plane orthogonal to~$\ba_i$
separates the edges into two parts: those oriented
into~$v_i$ from those oriented out of~$v_i$.  Thus,
there are exactly two inversions in this case.
\ $\sq$

\medskip

\subsection*{Final remarks}  \,
Let us note that Dehn's theorem and Cauchy's corollary hold
in the generality of all convex polytopes.  The  reduction
to the simplicial case is straightforward: triangulate the
surface of~$P$ by adding diagonals of the faces.  We skip
the details.  We should mention also that Dehn's original
proof is also based on a graph-theoretic argument, very much
different from the one presented here.

The proof we presented above easily splits into
two parts: global and local, not unlike the original
proof of the Cauchy theorem (see e.g.~\cite{Ber,C3}).
The local part (Lemma~2), while different, has the same
flavor as the the sign changes lemma in Cauchy's proof.
It is even more similar to the local part of other proofs
of Dehn's theorem (see~\cite{C3,Pak,Whi}).  The global
part is a graph-theoretic argument somewhat similar in
style and complexity to  the double counting
argument used in the proof of the Cauchy theorem, with
both arguments based on Euler's formula.

Finally, let us mention that the result of Trushkina is
a bit stronger as she allows the base vertices not to
be on the same face.  When applied to the infinitesimal
rigidity of polytopes, this gives an equivalent result,
but the graph-theoretic claim becomes harder to prove.
Similarly, Schramm proves a much stronger result about
certain delicate graph labeling~\cite{Sch}.
Interestingly, at some point he employs a counting angles
argument, essentially reproving Euler's formula.

\smallskip

\subsection*{Acknowledgements}  \, I am grateful to
Bob Connelly and Oded Schramm for their comments
on article~\cite{Sch}, and to the NSF for the
financial support.  I am also indebted to Victor
Alexandrov for his interest in the paper and
remarks on the presentation.  Last but not least,
I would like to thank Susanna Yusufova for her
love and support.

\medskip



\vskip.3cm

\textbf{\textsl{Igor Pak}}

\textsl{M.I.T.} \, \textsl{Department of Mathematics}

\textsl{Cambridge, MA 02139} \,

\textsl{U.S.A.}

\texttt{pak@math.mit.edu}


\begin{thebibliography}{12}\label{refpage}

\bibitem{Al-book}
A.~D. Aleksandrov, \emph{Vypuklye mnogogranniki {\rm (in
  Russian)}}, {M.: Gostekhizdat}, 1950; English translation:
\emph{Convex polyhedra}, Springer, Berlin, 2005.

\bibitem{Al-inverse} V.~A. Alexandrov,
Inverse function theorems and their applications to the theory
of polyhedra (in Russian), to appear in \emph{Reviews in
Mathematics and Mathematical Physics} (2005), 125 pp.;
available electronically at \ts {\tt http://math.nsc.ru/\~\/vaalex}

\bibitem{Ber} M. Berger, \emph{G\'{e}om\'{e}trie}, Vol. 1-5,
(in French), Nathan, CEDIC, Paris, 1977.

\bibitem{C3}
R. Connelly, Rigidity, in {\it Handbook of Convex Geometry},
vol.~A, 223--271, North-Holland, Amsterdam, 1993.

\bibitem{Dehn}
M. Dehn, \"{U}ber die Starreit konvexer Polyeder,
\emph{Math. Ann.}~\textbf{77} (1916), 466--473;
available electronically at \ts
{\tt http://dz-srv1.sub.uni-goettingen.de/cache/toc/D37460.html}

\bibitem{FP} M. Fedorchuk and I. Pak,
Rigidity and polynomial invariants of convex polytopes,
\emph{Duke Math.~J.} \textbf{129} (2005), 371–-404; available
electronically at \ts {\tt http://www-math.mit.edu/\~\/pak}

\bibitem{G}
H. Gluck, Almost all simply connected closed surfaces are rigid,
in {\it Lecture Notes in Math.}~\textbf{438}, 225--239, Springer,
Berlin, 1975.


\bibitem{Pak}
I. Pak, \emph{Lectures on Combinatorial Geometry and
Convex Polytopes}, monograph in preparation.


\bibitem{Roth}
B. Roth,
Rigid and flexible frameworks,
\emph{Amer. Math. Monthly}~\textbf{88} (1981), no. 1, 6--21.

\bibitem{Sch} O. Schramm, How to cage an egg,
\emph{Invent. Math.}~\textbf{107} (1992), 543--560;
available electronically at \ts
{\tt http://dz-srv1.sub.uni-goettingen.de/cache/toc/D183703.html}

\bibitem{Tru}
V.~I. Trushkina,
A coloring theorem and the rigidity of a convex polyhedron
(in Russian),
\emph{Ukrainian Geometric Sbornik}~No.~\textbf{24} (1981),
116--122.

\bibitem{Tru-next}
V.~I. Trushkina,
Method of 3-colouring of graphs,
\emph{Sib. Math. J.}~\textbf{28} (1987), 331--342.


\bibitem{Whi}
W. Whiteley,
Infinitesimally rigid polyhedra. I. Statics of frameworks,
\emph{Trans. AMS}~\textbf{285} (1984), 431--465.

\bibitem{Whi-survey}
W. Whiteley, Rigidity and scene analysis,
in \emph{Handbook of discrete and computational geometry}, 893--916,
CRC Press, Boca Raton, FL, 1997.


\end{thebibliography}
\end{document}